# OPTIMAL CONTROL OF SYSTEMS ENGINEERING. DEVELOPMENT OF A GENERAL STRUCTURE OF THE TECHNOLOGICAL CONVERSION SUBSYSTEM (PART 2)

I. Lutsenko
PhD, Professor
Department of Electronic Devices
Kremenchuk Mykhailo Ostrohradshyi
National University
Pervomaiskaya str., 20,
Kremenchuk, Ukraine, 39600
E-mail: delo-do@i.ua

*Запропоновано архітектуру технологічної підсистеми, особливістю якої є наявність максимальної кількості ступенів свободи. Це забезпечує можливість незалежної зміни всіх важливих параметрів технологічного процесу. Вбудовані можливості оцінки економічних параметрів технологічної операції вирішують завдання оптимізації технологічного процесу за критерієм ефективності використання ресурсів. Всі рішення перевірені з використанням спеціального системного конструктора EFFLI*

*Ключові слова: технологічна підсистема, структура оптимальної технологічної підсистеми, оптимізація технологічного процесу*

*Предложена архитектура технологической подсистемы, особенностью которой является наличие максимального количества степеней свободы. Это обеспечивает возможность независимого изменения всех значимых параметров технологического процесса. Встроенные возможности оценки экономических параметров технологической операции решают задачу оптимизации технологического процесса по критерию эффективности использования ресурсов. Все решения проверены с использованием специального системного конструктора EFFLI*

*Ключевые слова: технологическая подсистема, структура оптимальной технологической подсистемы, оптимизация технологического процесса*

## 1. Introduction

In [1], the basic structure of the sTransA conversion system with a batch feed of products on an example of a liquid heating system was developed. What is the feature of this structure?

Firstly, the conversion system solves **a simple task** – conversion of input raw products into the output product. In the demo, this is the conversion of cold liquid into a heated liquid.

It is necessary to pay attention to the phrase **"output product"**. This is the output product, not the **finished product** since the conversion system is only responsible for the qualitative indicator of the conversion product. Responsibility for the possibility of selecting the finished product with the required quantitative parameters by the consumer lies with the buffering system of the sSepA1controlled system.

Secondly, the conversion system implements the principle of **batch conversion** of raw products. At the output of the conversion system, a continuous finished product can not be obtained in principle. Such a product can be obtained at the output of the buffering system.

Third, principle of feed of volumes of raw products, independent of the energy product feed rate is realized in the conversion system. And, therefore, full implementation of optimal control technologies is potentially possible in the conversion systems, formed according to this principle. However, only potentially.

For the practical implementation of the optimal control principles, the basic structure of the conversion system should be supplemented by the necessary structures, which provide the delivery of information signals about the size of the resulting target product and the value of the optimization criterion for the technological operation (TO) performed.

The concept of "optimal" (the best), is very often, and it's an understatement, treated by researchers fairly arbitrary. Some believe that the best control is the control at minimum cost while others regard control with a maximum movement speed of the controlled object as the best, many of them minimize power consumption and so on.

In fact, this concept does not have "degrees of freedom". Control, which ensures the generation of maximum value added, which leads to maximum profits of the enterprise as a whole, if the optimization, by the efficiency criterion, covers all the systems of such enterprise is optimal.

The paper aims at developing a cybernetic structure of the technological subsystem, which is able to provide such optimization.

## 2. Analysis of the literature data and problem statement

The principle of determining the basic indicators of the technological process (TP) according to the results of the control of its parameters in an attempt to implement the optimal control principles has been always used, and began to be actively discussed in the special literature since the mid 19th century. However, the inadequate model of the technological conversion subsystem (TCS) (excluding the wear process) leads to false conclusions. For example, such statement is quite typical and still has not lost its "relevance": "Therefore, the operation of the circuit of enrichment devices can be optimized without the direct accounting of economic indicators. It is enough to provide a maximum output of the concentrate with regard to the restrictions, applied to the object that will provide the maximum profit" [2].

This approach in Soviet scientific schools and other has led to the development of entire areas, where problems of assessing technological processes were reduced to estimating the speed (performance) of the TCS [3, 4]. Accordingly, this was reflected in the model of the technological subsystem, the modes of which are not directly related to economic indicators.

Development of computer technology has led to the appearance of many publications, in which the optimization problem was formulated as the process of determining the minimum cost [5]. But the cost can be determined both by a separate technological process, and in relation to the larger structures of the enterprise, which seems easier. As can be seen from the description of the MRP II standard, standardization efforts are more focused on the structure of the enterprise as a whole, than the systems [6] it is composed of.

Failures in attempts to use the control at minimum cost have led to the fact that the TCS models, as a rule, are developed based on the physicochemical aspects of the processes [7].

An important factor that hinders gathering statistical data for making managerial decisions directly within the TCS is measurement errors of quantitative parameters of conversion products. This is a strong reason to use the optimization criterion, since the use of direct methods for determining the optimum requires by orders of magnitude more experiments with the need to process the "cloud" of measurements [8].

Problems, associated with the lack of a cybernetic structure of the TCS lead to attempts to create genetic algorithms [9]. It is expected that such an algorithm, in the course of evolution, generates the optimal TCS structure.

## 3. Goal and objectives

The goal of the work is to develop the TCS architecture, which provides the formation of not only the output consumer product, but also information products, displaying the absolute and the relative indicators of the extent of achieving the control objectives.

Tasks that need to be solved to achieve this goal, are:
– formation of basic subsystems of the conversion system;
– development of mechanisms for the superstructure of the control subsystem that allows to test the capabilities of the technological conversion subsystem of products and input mechanisms;
– formation of the TCS that provides possibilities for extreme (optimal) control.

## 4. Formation of the basic subsystems of the conversion system

Synthesis of the basic internal structure of the conversion system [1] allows to proceed to the formation of the subsystems of the system. In the current scientific and engineering works, related to the study of the internal structure of the cybernetic system, there is no description of the reasonable concept of reference of this or that object of the system to its subsystems.

Experience of practical integration of internal objects of the investigated systems shows that integration of the internal system objects according to the principle of the intensity of the interaction of these objects is the most natural. That is, the structures, among which there is the intensive data exchange should not be in different subsystems.

This is quite natural, since experience shows that it is more effectively to place employees, intensively interacting in working moments in the same office.

In the formed system [1], synthesized mechanisms ensure the functioning of the technological process (TP) and the control process. This means that it is necessary to form two subsystems: technological subsystem and control subsystem within the conversion system.

Based on the object-oriented principle, the mTmprA1 mechanism and mechanisms that continually monitor the process parameters (mFinA1, mCmpA1, mFinA2), should be attributed to technological subsystem. All the rest - to the control subsystem.

In this case, one unit pulse signal is transmitted for one functioning cycle of the system via RTF-RTF, RED-RED and PTF-PTF channels (Fig. 1).

Signal that is transmitted via RTF-RTF channel, informs the control subsystem of the completion of receipt of the raw product in the mTmprA1 object buffering mechanism.

Signal that is transmitted via RED-RED channel, notifies the control subsystem that the quality of the conversion product has reached the set value.

Finally, the signal transmitted via PTF-PTF channel, indicates that the finished product is transferred to the buffering system.

By combining the mTmprA1, mFinA1, mCmpA1 and mFinA2 mechanisms into the sbTransA technological subsystem, and mPassA1, mPassB1 and mGstB1 mechanisms into the sbContrA1 control subsystem, we get the opportunity to present the heating system in the form of two subsystems – the technological subsystem and control subsystem (Fig. 2).

Since the control unit is not an integral object, the term **control system** is apparently incorrect. Using it can be justified only by the historical practice.

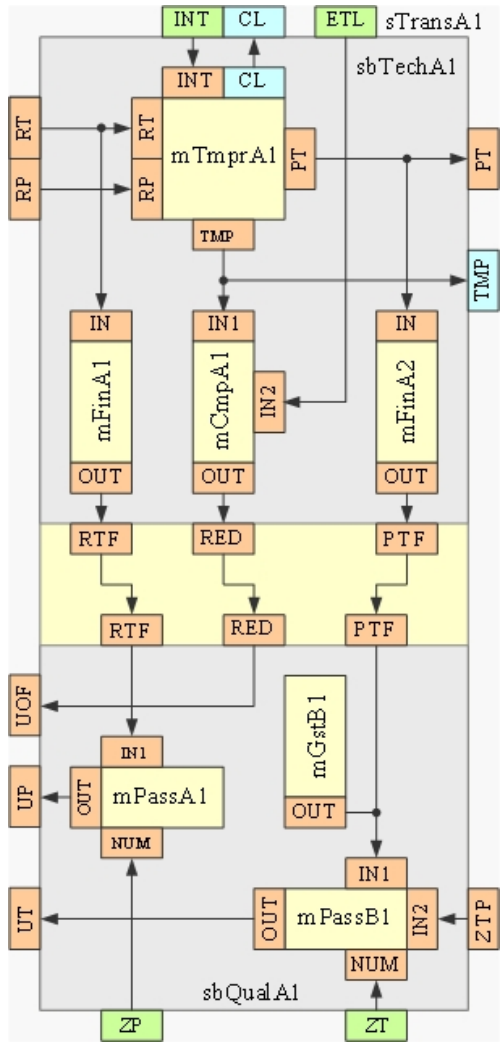

Fig. 1. The result of integration of the internal structure of the system under study according to the object-oriented principle

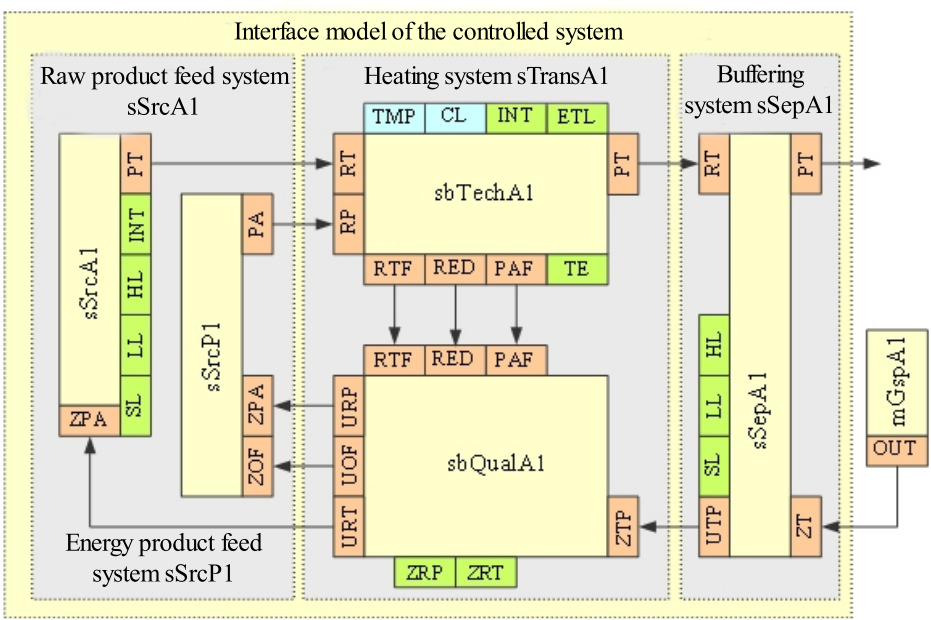

Fig. 2. The interface model of the controlled system, conversion system of which is represented in the form of liquid heating subsystem and control subsystem

The architecture of the sTransA basic system has great potential with respect to the possibility of setting various control modes. Thus, we can control the energy product feed rate, change the raw product (cold liquid) feed amount. The quality of the output product will be the same, within a certain control range.

However, it will be impossible to implement the potential of the system in relation to the formation of the target product if you can not control the time of TO, quantitative parameters of the input and output products of TO, integrate them, reduce to comparable values and use that information to form the optimal control criterion.

Advanced controlled systems should not only form the consumer product with the required qualitative and quantitative parameters, but also ensure the maximization of the target product (value added in economic systems or value added in cybernetic systems).

## 5. Investigation of control modes

### 5. 1. Development of additional objects

Let us consider the change of the energy product consumption in the system in the control process and time of technological operation.

To conduct the study, additional **EFFLI** objects: linear scanner of the control range, integrator with reset, timer, report generator were created.

### 5. 1. 1. The mScanB linear scanner of the control range

The mScanB linear scanner of the control range is designed to automatically change the control signals with a given direction and a predetermined step in a certain range. The interface model of the control scanner is shown in Fig. 3.

Assignment of the mScanB port sections is given in Table 1.

Table 1
Assignment of mScanB port sections

| Port assignment | Notation |
|---|---|
| Control signal minimum | MIN |
| Control signal maximum | MAX |
| Control change step | STP |
| Control range scanning direction | DIR |
| System stop when testing the control range | STS |
| Strobe signal | STR |
| Control signal | OUT |
| The signal of reaching the control range boundary | RPT |

The first five parameters are set before the start of the control process.

If a zero value is set in the DIR section, the control changes from a smaller value (set in the MIN section) to the larger value (set in the MAX section). Unit level in the DIR section leads to the control change from a larger value to smaller.

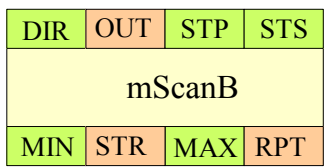

Fig. 3. The interface model of the mScanB linear scanner of the control range in the form of EFFLI object

Control step is set in the STP section.

Setting the unit-level signal in the STS section provides the termination of the operation of the controlled system after reaching the control boundary.

Control signal change in the OUT section occurs when the unit-level strobe signal enters the STR port section.

### 5. 1. 2. The mTimerA timer

Time of the technological operation (TO) is a necessary parameter for solving the optimal control problem and mTimerA timer was created for its determination. The interface model of the timer is shown in Fig. 4.

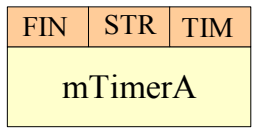

Fig. 4. The interface model of the mTimerA timer in the form of EFFLI object

Assignment of the mTimerA timer port sections is shown in Table 2.

Table 2

Assignment of the mTimerA timer port sections

| Port assignment | Notation |
| --- | --- |
| Timer start | STR |
| Timer stop | FIN |
| Measured time interval | TIM |

### 5. 1. 3. The mIntA integrator

Determining the energy product consumption requires a mechanism that ensures the conversion of the energy flow into the energy consumption volume. For these purposes, the mIntA integrator with reset was developed. The interface model of the integrator is shown in Fig. 5.

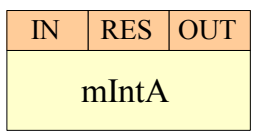

Fig. 5. The interface model of the integrator mIntA in the form of EFFLI object

Assignment of the mIntA integrator port sections is given in Table 3.

Table 3

Assignment of the mIntA integrator port sections

| Port assignment | Notation |
| --- | --- |
| Input signal | IN |
| Output signal | OUT |
| Integrator reset signal | RES |

### 5. 1. 4. The mBaseA report generator

The data needed for the analysis, obtained during the operation of the controlled system should be presented in a readable form. For these purposes, the mBaseA report generator was developed. The interface model of the report generator is shown in Fig. 6.

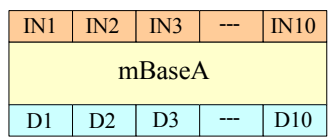

Fig. 6. The interface model of the mBaseA report generator in the form of EFFLI object

Assignment of the mBaseA report generator port sections is given in Table 4.

Table 4

Assignment of the mBaseA report generator port sections

| Port assignment | Notation |
| --- | --- |
| Record number | NUM |
| Strobe signal | STR |
| Input signal 1 | IN1 |
| Input signal 2 | IN2 |
| ----------------------- | ---- |
| Input signal 10 | IN10 |
| Output signal 1 | OUT1 |
| Output signal 2 | OUT2 |
| ----------------------- | ------- |
| Output signal 10 | OUT10 |

During the next operation, signals, entering the input section are saved, and upon the completion of the operation (upon arrival of the strobe signal) are transferred to the output, after which, if required, displayed on the report sheet. The NUM section generates the operation number.

### 5. 2. The model of the controlled system with the possibility to determine additional parameters

Assembled model (Fig. 7) in a software designer in the form of EFFLI objects allows to display the process of change of energy consumption volumes of operations and change in the time of TO in the process of the control change.

The results of the operation of the controlled system in the form of diagrams of the change of energy consumption and time of TO are shown in Fig. 8, 9.

As can be seen from the timing diagrams, the time of TO and energy costs for its implementation decrease with an increase in the energy product feed rate.

Based on the analysis of timing diagrams, it seems that the higher the energy product feed rate, the better economic indicators of TO.

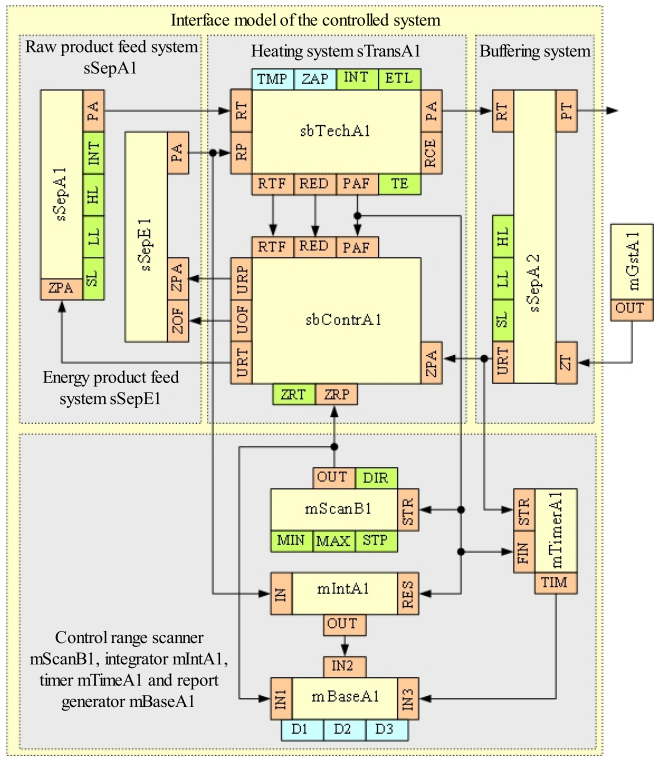

Fig. 7. The architecture of the controlled system for the study of the effect of control modes on the energy product consumption and operation time

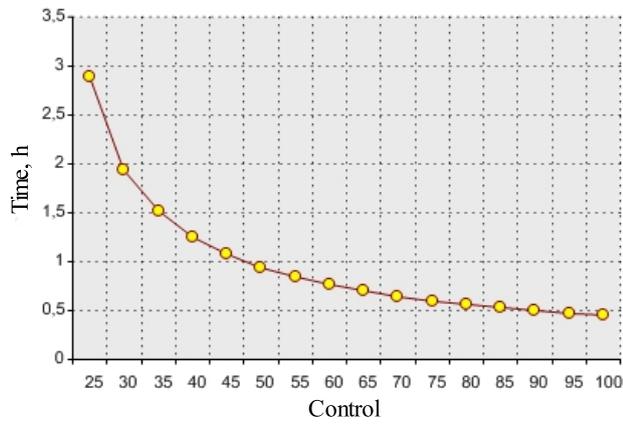

Fig. 8. Change of the operation time depending on the control

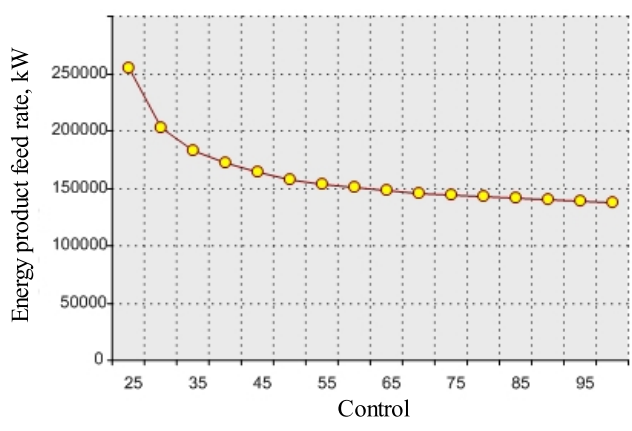

Fig. 9. Change of energy consumption of operations depending on the control

### 5.3. Development of the structure of the technological subsystem with the optimization possibility

The physical mechanism of the above processes is clear. The higher the energy product feed rate, the faster the liquid is heated, and the decrease in the time of TO naturally leads to the reduction of heat loss during the heating process.

However, in this case, the most favorable mode is the mode with maximum energy product feed rate.

All this means only that the model of the TP of heating is not perfect. The wear of the heating mechanism [10] increases disproportionately with the increase in the heating process speed. According to this approach, in the developed model of the TCS, the service life of the object depending on the energy product feed rate is described by the expression

$$T = T_n k^{-\alpha},$$

where $T$ is the service life of the object (mean time between failures) when operating in a given mode; $T_n$ is the service life of the object when operating in the nominal mode; $k$ is the ratio of the actual energy product feed rate to the feed rate in the nominal mode; $\alpha$ is an indicator, associated with the used conversion technology of a special product.

To form the current value of the heating mechanism wear rate, the mWearA mechanism was developed (Fig.10)

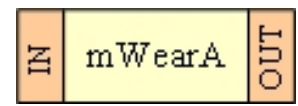

Fig. 10. The mechanism of the formation of the heater wear signal

When delivering data of current energy product rate to the IN input, mWearA generates signal of the current wear rate of the heating mechanism on the OUT output. The higher the energy product feed rate, the disproportionately higher the heating mechanism wear. This fact must be taken into account, because otherwise we will not see a full, natural picture of the physics of the process, and, accordingly, will not be able to determine the optimal operation mode of the TP.

However, in order to purposefully change the operation mode of the controlled system, it is necessary to compare not the actual values of energy consumption and wear, but expenses and income from the TO with respect to the operation time.

The process of comparison is fairly simple. It is necessary to determine the amounts of raw and energy product, as well as the wear of the technological mechanism for the operation time. By determining their cost estimate, we obtain the cost of the technological operation. These costs are equal to the sum of costs for each type of input products. Also, the cost estimate of the output product of the operation can be obtained in a similar way.

For the technical implementation of this ideology, it is necessary to create several more mechanisms: multiplier (mMultA), summator (mSumA) and TO identifier (mIdsA) (Fig. 11).

Assignment of the multiplier and summator sections is common and does need the explanation.

Assignment of the TO identifier port sections is shown in Table 5.

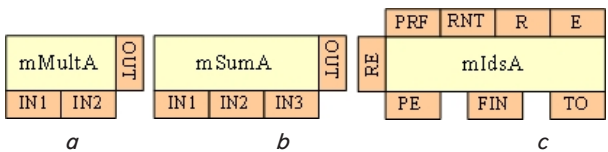

Fig. 11. The interface model of mechanisms *a* — mMultA, *b* — mSumA, *c* — mIdsA in the form of EFFLI objects

Table 5

Assignment of the mIdsA TO identifier port sections

| Port assignment | Notation |
|---|---|
| Cost estimate of input products of the operation | RE |
| Cost estimate of output products of the operation | PE |
| The time of the beginning of the calculation of indicators | FIN |
| Technological operation time | TO |
| Value added of the operation | PRF |
| Profitability of the operation | RNT |
| Resource intensity of the operation | R |
| Efficiency of the operation | E |

The interface model of connection of computational mechanisms to identify TO is shown in Fig. 12.

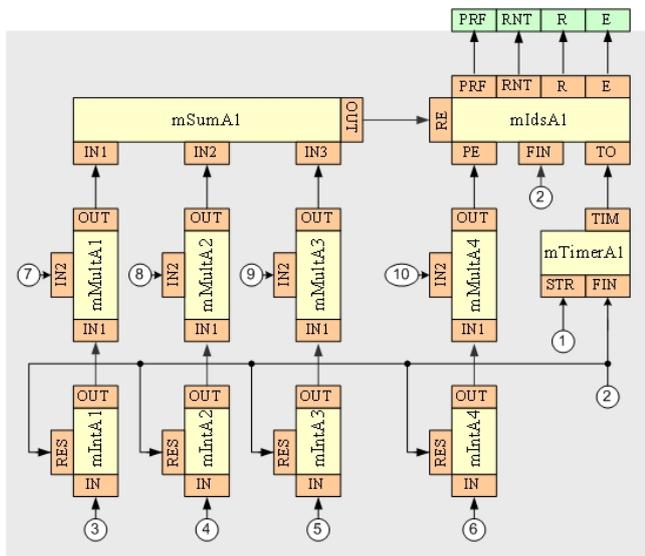

Fig. 12. Connection diagram of computational mechanisms: 1 — start of the raw product feed; 2 — completion of the output product release; 3 — signal of the energy product feed rate; 4 — signal of the raw product feed rate; 5 — signal of the wear rate; 6 — signal of the output product feed rate; 7 — cost estimate of the unit of the energy product; 8 — cost estimate of the unit of the raw product; 9 — cost estimate of the unit of the wear; 10 — cost estimate of the unit of the output product

Since, in this case, from the operation to operation, only energy product consumption changes, the cost estimate of the output product of the operation is unchanged.

In the technological subsystem, upon the completion of the TO, indicators, relating to the target product (value added) are determined at a certain time interval based on three key parameters.

Conversion subsystem, which provides the implementation of all the necessary operations, is shown in Fig. 13.

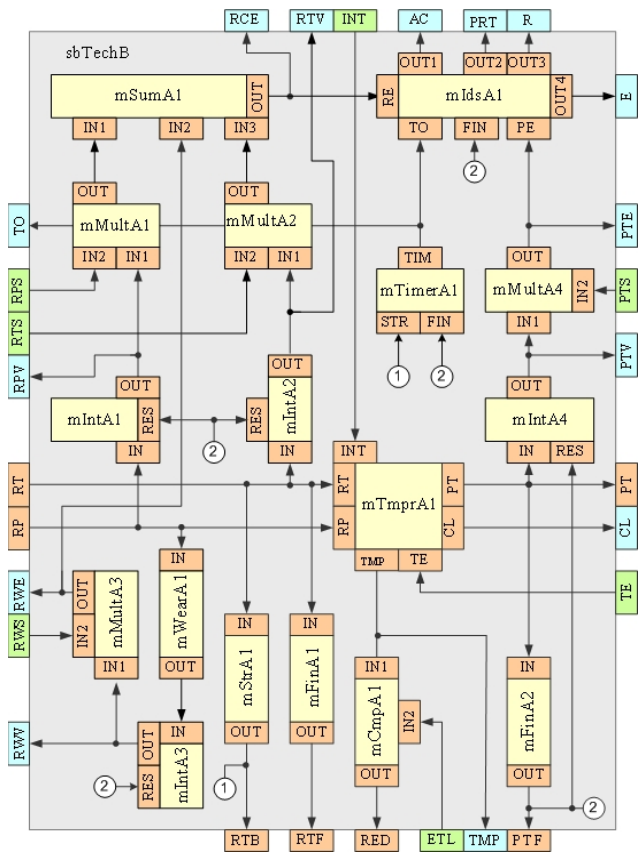

Fig. 13. The structure of the technological subsystem in the form of EFFLI objects

Assignment of the sbTechB technological subsystem port sections is given in Table 6.

To test the TCS settings, dependencies of the heating mechanism wear and energy consumption on the control were removed. The nature of the dependencies shows that, in the control process, the cost function extremum will be within the selected control range (Fig. 14).

Fig. 15 shows the timing diagrams of the change of the cost estimate of costs, cost estimate of the output product of operations and the time of operations depending on the control.

Changing cost estimates of raw, energy and output product on sheets of corresponding objects, displacement of the minimum cost and line of the cost estimate of the output product can be observed.

As will be shown below, determination of basic indicators $RE$, $PE$ and $T_{op}$ provides the possibility of obtaining all derivatives of economic indicators, needed for determining the value of the target product of the TO and solving the optimal control problem. In the R and E sections on the basis of the key indicators, "resource intensity" [11] and "efficiency" [12] indicators, which provide the optimization mode are formed.

Assembled model effli_model_2 can be downloaded here [13]. After the start, report generator forms the data, which were used for building the diagrams (Fig. 14, 15).

Table 6

Assignment of the sbTechB technological subsystem port sections

| Port assignment | Notation |
|---|---|
| Output product release intensity | INT |
| Ambient temperature | TE |
| Heating temperature setpoint | ETL |
| Cost of the unit of the raw product | RTS |
| Cost of the unit of the energy product | RPS |
| Cost of the unit of the equipment service life | RWS |
| Cost of the unit of the output product | PTS |
| Raw product | RT |
| Energy product | RP |
| Output product | PT |
| Start of receipt of the raw product | RTB |
| Completion of receipt of the raw product | RTF |
| Heating process completion signal | RED |
| Output product release completion | PTF |
| Operation time | TO |
| Wear amount | RWE |
| Amount by the output product | PTE |
| Cost amount | RCE |
| Value added of the operation | AC |
| Profitability of the operation | PRT |
| Efficiency of the operation | E |
| Sign of the presence of losses of start | RWM |
| Raw product consumption volume | RTV |
| Energy product consumption volume | RPV |
| Output product release volume | PTV |
| Equipment wear for the operation time | RWV |
| Current value of the heating temperature | TMP |
| Heating mechanism load level | CL |
| Resource intensity | R |

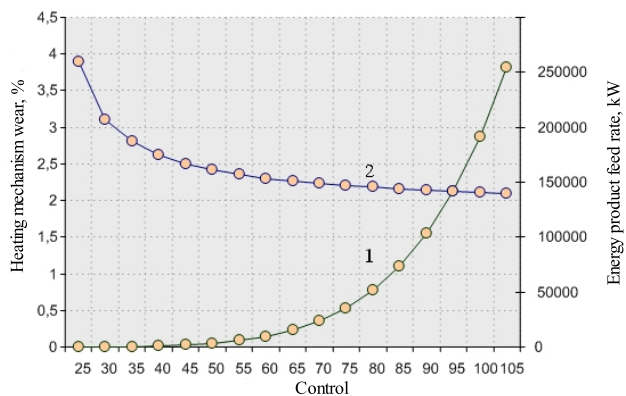

Fig. 14. Diagrams of the change of:
1 — energy costs and 2 — technological equipment wear

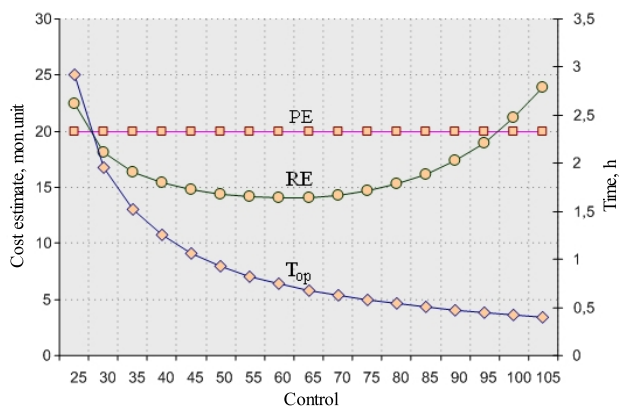

Fig. 15. Diagrams of the change of the cost estimate of costs (RE), cost estimate of the output product of operations (PE) and time of operations depending on the control ($T_{op}$)

## 7. Conclusions

In this part of the work, previously developed structure of the basic conversion system has been divided into two subsystems: the basic technological subsystem and control subsystem. Distribution of objects was carried out according to the object-oriented principle.

Research of energy consumption of technological operations depending on the control has shown the need to finalize the internal structure of the technological subsystem. The availability of accounting the technological equipment wear is of fundamental importance since it is often impossible to determine the extremum of losses or determine it reliably in the other cases.

In addition to a consumer product, target product (value added) is always the output product of the technological subsystem, the technological subsystem should be able to provide information on the value of the product and the value of the efficiency criterion. Without this information, it is impossible to ensure the implementation of optimal control.

The maximum number of degrees of freedom together with the ability to assess technological operations opens up the practical possibility of implementing new, previously unavailable in terms of effectiveness technologies in production.

# ПОСТРОЕНИЕ СЛЕДЯЩИХ ИНВАРИАНТНЫХ СИСТЕМ НА ОСНОВЕ ЭКВИВАЛЕНТНОГО РОБАСТНОГО УПРАВЛЕНИЯ

*Пропонується підхід до побудови інваріантних відносно зовнішніх збурень систем управління. На відміну від класичної схеми тут не потрібно вимірювати обурення. Можливість безмежного збільшення коефіцієнта посилення робастного регулятора без втрати стійкості дозволяє зменшити вплив збурень до як завгодно малої величини. Цим забезпечується висока точність стеження еталонної траєкторії і швидкодія для широкого класу збурень*

*Ключові слова: система управління, інваріантна система, невизначеність, функція Ляпунова, робастний регулятор, обурення*

*Предлагается подход к построению инвариантных относительно внешним возмущениям систем управления. В отличие от классической схемы здесь не требуется измерять возмущения. Возможность беспредельного увеличения коэффициента усиления робастного регулятора без потери устойчивости позволяет уменьшить влияние возмущений до сколь угодно малой величины. Этим обеспечивается высокая точность слежения эталонной траектории и быстродействие для широкого класса возмущений*

*Ключевые слова: система управления, инвариантная система, неопределенность, функция Ляпунова, робастный регулятор, возмущение*

**Г. А. Рустамов**
Доктор технических наук, профессор
Кафедра «Автоматика и управление»
Азербайджанский Технический Университет
пр. Г. Джавида, 25,
г. Баку, Азербайджан, Az1073
E-mail: gazanfar.rustamov@gmail.com

## 1. Введение

Одним из эффективных и простых методов повышения точности систем автоматического регулирования является принцип инвариантности или принцип управления по возмущению [1–5]. При этом основной задачей является обеспечение инвариантности относительно к внешним возмущениям, действующих на объект.

Следует отличать астатические системы от инвариантных. В астатических системах возмущение не измеряется и полная компенсация статической ошибки $\Delta_s$ за счёт обратной связи достигается в пределе: $\Delta s = \lim e(t)$ при $t \to \infty$.

В инвариантных системах компенсация возмущения происходит почти мгновенно. Однако построение хотя приближенно инвариантной системы, в которой удовлетворяется лишь нестрогие (ослабленные) условия инвариантности, за счёт обратной связи наталкивается на принципиальные трудности. В классических задачах эта проблема решена путем измерения возмущения и воздействия на тот же объект через *инверсный регулятор* (компенсатор), находящийся во втором канале вне обратной связи.

Как и во многих областях, комбинированный подход несколько улучшает положение: за счёт обратной связи компенсируется неконтролируемые возмущения, а подавление доминирующих измеряемых возмущений возлагается на компенсатор [3].

На практике «измерение» в зависимости от физической природы возмущения может оказаться проблематичным из-за следующих причин: отсутствие методики измерения; сложность и слабое быстродействие измерительной системы, а также её габариты и вес;